\documentclass[12pt]{article}
\usepackage[utf8]{inputenc} %for Windows
\usepackage[german, english]{babel}
\usepackage[tbtags]{amsmath}
\usepackage{amsfonts,amssymb,mathrsfs,amscd,comment}
\usepackage{tikz}
\usepackage{amsthm}
\usepackage[a4paper, total={6in, 9in}]{geometry}

\DeclareMathOperator{\mmod}{mod}

\DeclareMathOperator{\Pol}{Rat}
\DeclareMathOperator{\Log}{Log}

\newtheorem{lemma}{Lemma}
\newtheorem*{lemma*}{Lemma}
\newtheorem{theorem}{Theorem}
\newtheorem*{theorem*}{Theorem}

\newtheorem{proposition}{Proposition}
\newtheorem{remark}{Remark}
\newtheorem*{remark*}{Remark}

\newcommand*{\bfrac}[2]{\genfrac{}{}{0pt}{}{#1}{#2}}
{}{}{}

\global\long\def\M{\mathcal{M}}
\global\long\def\N{\mathcal{N}}
\global\long\def\S{\mathcal{S}}
\global\long\def\Rc{\mathcal{R}}

\title{A variant of the Bombieri-Vinogradov theorem with explicit constants}
\date{}
\author{Alisa Sedunova}

\begin{document}
\maketitle

\begin{abstract} 
In this paper we improve the result of \cite{Akbary2015} with getting $(\log x)^{\frac{7}{2}}$ instead of $(\log x)^{\frac{9}{2}}$. In particular we obtain a better version of Vaughan's inequality by applying the explicit variant of an inequality connected to the Möbius function from \cite{Helfgott2014a}.
\end{abstract}

\section{Introduction}
For integer number $a$ and $q \geq 1$, let 
$$
%\psi(x)=\sum_{\substack{n\leq x }} \Lambda(n) \qquad \text{ and } \qquad 
\psi(x; q, a)=\sum_{\substack{n\leq x \\ {n\equiv a\bmod{q}}}} \Lambda(n),
$$
where $\Lambda(n)$ is the von Mangoldt function. 
The Bombieri-Vinogradov theorem is an estimate for the error terms in the prime number theorem for arithmetic progressions averaged over all $q \leq x^{1/2}$. 
\begin{theorem*}{\bf (Bombieri-Vinogradov)}
Let $A$ be a given positive number and $Q\leq x^{1/2} (\log{x})^{-B}$ where
$B=B(A)$, then
$$\sum_{q \leq Q} \max_{2 \leq y \leq x} {\max_{\bfrac{a}{(a,q)=1}}} \left|\psi(y,q,a) - \frac{y}{\phi(q)} \right| \ll_A \frac{x}{(\log{x})^A}.$$
\end{theorem*}
The implied constant in this theorem is not effective, since we have to take care of characters, associated with those $q$ that have small prime factors.
The main result of this paper is
\begin{theorem} {\bf (Bombieri-Vinogradov theorem with explicit constants)} \label{bvef}
Let $x \geq 4$, $1 \leq Q_1 \leq Q \leq x^{\frac{1}{2}}$. Let also $l(q)$ denote the least prime divisor of $q$. Define $F(x,Q,Q_1)$ by 
\begin{equation*}
\begin{split}
& F(x,Q,Q_1)=\frac{14x}{Q_1}+4x^{\frac{1}{2}}Q+15x^{\frac{2}{3}}Q^{\frac{1}{2}}+4x^{\frac{5}{6}}\log \frac{Q}{Q_1}.
\end{split}
\end{equation*}
Then 
$$
\sum_{\bfrac{q\leq Q }{ l(q)>Q_1}} \max_{2 \leq y \leq x} \max_{\bfrac{a}{(a,q)=1}} \left| \psi(y;q,a) - \frac{\psi(y)}{\phi(q)} \right| 
< c_1 F(x,Q,Q_1)(\log x)^{\frac{7}{2}}, $$
where 
\begin{equation*}
\begin{split}
& c_1= \frac{5}{4} E_0 c_0+1 = 42.140461 \ldots,\\
& E_0 = \prod_{p} \left(1 + \frac{1}{p(p-1)}\right) = 1.943596\ldots,\\
& c_0 = (2A_0)^{\frac{1}{2}}\frac{2^5}{3^{\frac{3}{2}}\pi (\log 2)^2} \left(2+\frac{\log(\log2)}{\log\frac{4}{3}}\right) = 16.93375 \ldots,\\
& A_0 = \max_{x >0} \left(\frac{\psi(x)}{x}\right) = \frac{\psi(113)}{113} = 1.03883\ldots.\\
\end{split}
\end{equation*}
\end{theorem}
Previously the best result obtained by these methods in the literature is due to Akbary, Hambrook (see \cite[Theorem 1.3]{Akbary2015}), where they proved that under assumptions of Theorem \ref{bvef} we have. 
$$
\sum_{\bfrac{q\leq Q }{ l(q)>Q_1}} \max_{2 \leq y \leq x} \max_{\bfrac{a}{(a,q)=1}} \left| \psi(y;q,a) - \frac{\psi(y)}{\phi(q)} \right| 
< c_1 F(x,Q,Q_1)(\log x)^{\frac{9}{2}}, $$
where $F(x,Q,Q_1)$ is defined by 
\begin{equation*}
\begin{split}
& F(x,Q,Q_1)=\frac{4x}{Q_1}+4x^{\frac{1}{2}}Q+18x^{\frac{2}{3}}Q^{\frac{1}{2}}+5x^{\frac{5}{6}}\log \frac{eQ}{Q_1}.
\end{split}
\end{equation*}
Here we reduce this power to $(\log x)^\frac{7}{2}$ by applying
an explicit version for an upper bound for 
$$b_k=\sum_{\bfrac{d \leq V}{d|k}}{\mu(d)},$$
where $\mu(d)$ is Mobius function, $V$ is a given number. This version can be found in \cite{Helfgott2014a}, namely  we have
\begin{lemma}{\bf (Helfgott, \cite[(6.9), (6.10)]{Helfgott2014a})} \label{hh}
For $V$ large enough we have
$$\sum_{k \leq Y}|b_k|^2 = Y(L+O^*(C))+O^*(V^2), \text{   where   }C=0.000023,L=0.440729$$
and $O^*(x)$ means that it is less in absolute value than $x$.
\end{lemma}
This Lemma is a variant of the sum considered in \cite{DIT}, where it is shown that
$$\sum_{d_1,d_2 \leq Y} \frac{\mu(d_1)\mu(d_2)}{\gcd(d_1,d_2)}$$
tends to a positive constant as $Y \to \infty$. It is also suggested without proving that $L$ can be about $0.440729$.

Notice, that by sharpening the inequality in Lemma \ref{hh} we will not be able to reduce the power of $\log x$, since the upper bound is optimal there, so by these methods the power $\frac{7}{2}$ is the best possible. Going further seems to be a hard problem which involves among simpler things a very careful analysis of the logarithmic mean of Möbius function twisted by a Dirichlet character.

\begin{remark*}
Let $Q=\frac{x^{\frac{1}{2}}}{(\log x)^B}$, where $B>\frac{7}{2}$. Then Theorem \ref{bvef} gives us the following bound 
\begin{equation*}
\begin{split}
&\sum_{\bfrac{q\leq Q }{ l(q)>Q_1}} \max_{2 \leq y \leq x} \max_{\bfrac{a}{(a,q)=1}} \left| \psi(y;q,a) - \frac{\psi(y)}{\phi(q)} \right| < c_1 \left(\frac{14x}{Q_1} (\log x)^\frac{7}{2}+19x(\log x)^{\frac{7}{2}-B}\right).
\end{split}
\end{equation*}
\end{remark*}
\begin{remark*}
It would be very good for applications to get $(\log x)^2$ in Theorem \ref{bvef}, however it seems impossible to get by present methods.
\end{remark*}

\begin{remark} \label{pi}
Define
$$
\pi(x) = \sum_{p \leq x} 1 \qquad \text{ and } \qquad \pi(x;q,a) = \sum_{\bfrac{p \leq x}{ p \equiv a \bmod{q} }} 1.
$$
Then Theorem \ref{bvef} under the same assumptions can be also formulated for $\pi(x)$, $\pi(x;q,a)$:
$$
\sum_{\bfrac{q\leq Q }{ l(q)>Q_1}} \max_{2 \leq y \leq x} \max_{\bfrac{a}{(a,q)=1}} \left| \pi(y;q,a) - \frac{\pi(y)}{\phi(q)} \right| 
< c_2 F(x,Q,Q_1)(\log x)^{\frac{7}{2}}, $$
where $c_2=1+\frac{2c_1}{\log 2}$.
\end{remark}
Proof of the remark is exactly the same as in \cite{Akbary2015}, we just have to change the power of $\log$.
%\begin{remark*}
%The condition $l(q)>Q_1$ in Theorem \ref{bvef} is taken to get rid of Siegel zeroes of the corresponding $L$-functions. It means that for $Q_1$ there exists just one primitive real character $\chi_{Q_1}$ such that $L(s, \chi_{Q_1})$ has the zero for $\Res s \geq 1 - \frac{c}{\log Q_1}$. By Siegel's theorem it follows that for $Q \geq Q_1$ the corresponding $L(s,\chi_Q)$ has zeros just with $\Res s \leq 1-\frac{c(\delta)}{Q^{\delta}}$ for $\delta >0$.
%\end{remark*}

The key tool for the proof of Theorem \ref{bvef} is Vaughan's identity, which we have to get in an explicit version for our goal.
Define 
$$\psi(y,\chi)= \sum_{n\leq y} \Lambda(n) \chi(n),$$
the twisted summatory function for the von Mangoldt function $\Lambda$ and a Dirichlet character $\chi$ modulo $q$.
One of two main results of this paper is
\begin{proposition}{\bf (Vaughan's inequality in an explicit form)} \label{vaug}
For $x \geq 4$
\begin{equation*}
\sum_{q \leq Q} \frac{q}{\phi(q)} \sum^*_{\chi (q)} \max_{y \leq x} |\psi(y,\chi)| <  c_0 (7x+2Q^2x^{\frac{1}{2}}+5Q^{\frac{3}{2}} x^{\frac{2}{3}}+4Qx^{\frac{5}{6}})(\log x)^{\frac{5}{2}},
\end{equation*}
where $Q$ is any positive real number and $\sum^{*}_{\chi (q)}$ means a sum over all primitive characters $\chi (\mmod q)$.
\end{proposition} 
The goal is to get an explicit version of $f(x,Q)$ by applying an improved version of P\'{o}lya-Vinogradov inequality (see \cite{Pomerance2010}), that will reduce the coefficients of $f(x,Q)$ and then we can apply Lemma \ref{hh}.

\section{Proof of Proposition 1}
Fix arbitrary real numbers $Q > 0$ and $x \geq 4$.
In this section, we shall establish Proposition \ref{vaug}, which is the main ingredient in the proof of Theorem \ref{bvef}.
Here we follow the ideas of \cite{Akbary2015} and applying the results from \cite{Helfgott2014a}. 
The main tool in the proof is the large sieve inequality (see, for example \cite[p.561]{Montgomery1978})
\begin{align} \label{sieve0}
\sum_{q \leq Q} \frac{q}{\phi(q)} {\sum^*_{\chi (q)}} \left| \sum_{m = m_0 + 1}^{m_0 + M} a_m \chi(m) \right|^2
\leq (M + Q^2) \sum_{m = m_0 + 1}^{m_0 + M} |a_m|^2,
\end{align}
from which it follows (see \cite[Lemma 6.1]{Akbary2015}) that

\begin{equation} \label{sieve}
\begin{split}
&\sum_{q \leq Q} \frac{q}{\phi(q)} {\sum^*_{\chi (q)}} \max_{y} \left| \sum_{m = m_0}^{M} \sum_{\bfrac{n = n_0}{mn \leq y}}^{N} a_m b_n \chi(mn) \right| \leq \\
& c_3 (M' + Q^2)^{\frac{1}{2}}(N' + Q^2)^{\frac{1}{2}} \left(\sum_{m = m_0}^{M} |a_m|^2\right)^{\frac{1}{2}}\left(\sum_{n = n_0}^{N} |b_n|^2\right)^{\frac{1}{2}} L(M,N),
\end{split}
\end{equation}
where 
%5\begin{equation} \label{c3}
%c_3 = \frac{2}{\pi} \left( \frac{2 + \log \frac{\log 2}{\log 4/3}}{\log 2}\right) = 2.64456\ldots,
%\end{equation}
$c_3=2.64...$, $L(M,N)=\log (2MN)$ and $M'=M-m_0+1$, $N'=N-n_0+1$ are the number of terms in the sums over $m$ and $n$ respectively.
Here the $a_m$, $b_n$ are arbitrary complex numbers.

\subsection{Sieving and Vaughan's identity}
We reduce to the case $2 \leq Q \leq x^{1/2}$.  If $Q < 1$, then the sum on the left-hand side of (\ref{vaug}) is empty and we are done.  Next, $1 \leq Q < 2$ then only the $q = 1$ term exists and we have
\begin{equation}
\sum_{q \leq Q} \frac{q}{\phi(q)} \sum^{*}_{\chi (q)} \max_{y \leq x} | \psi(y,\chi) | = \max_{y \leq x} \left| \sum_{n \leq y} \Lambda(n)  \right| = \psi(x) \leq A_0 x,
\end{equation}
which is better than the theorem. Finally, if $Q > x^{1/2}$, Theorem \ref{vaug} follows from (\ref{sieve}) with $M = m_0 = n_0 = 1$, $N = \lfloor x \rfloor$, $a_m = 1$, $b_n = \Lambda(n)$ by the estimate
$$
\sum_{n \leq x} \Lambda(n)^2 \leq \psi(x) \log x \leq A_0 x \log x.
$$

From now on we assume $2 \leq Q \leq x^{1/2}$.
Notice that the fact that we can restrict ourselves to the range $2 \leq Q \leq x^{1/2}$ allows us to apply Lemma \ref{hh}(otherwise it would make less sense, since the main term in Lemma \ref{hh} would be smaller than $O^*$-term).
As in \cite{Akbary2015} we will use Vaughan's identity (see also \cite{Vaughan1977})
$$
\Lambda(n) = \lambda_1(n) + \lambda_2(n) + \lambda_3(n) + \lambda_4(n),
$$
where 
\begin{equation*}
\begin{split}
&\lambda_{1}(n) =\begin{cases}
\Lambda(n), &\mbox{ if } n \leq U,\\
0, &\mbox{ if } n > U,
\end{cases}
\;\;\;\;\; \lambda_{2}(n) = \sum_{\bfrac{hd = n}{d \leq V}} \mu(d) \log h,\\
&\lambda_{3}(n) = - \sum_{\bfrac{mdr = n}{m \leq U, d \leq V}} \Lambda(m) \mu(d), \;\;\;
\lambda_{4}(n) = - \sum_{\bfrac{mk = n}{m > U,k > V}} \Lambda(m)  \sum_{\bfrac{d \mid k} {d \leq V}} \mu(d).\\
\end{split}
\end{equation*}

Assume $y \leq x$, $q \leq Q$, and $\chi$ is a character mod $q$.  We use the above decomposition to write
$$
\psi(y,\chi) = S_1 + S_2 + S_3 + S_4,
$$
where
$$
S_i = \sum_{n \leq y} \lambda_i(n) \chi(n).
$$
Let $U$, $V$ be non-negative functions of $x$ and $Q$ to be set later and denote the contributions to our main sum by
$$\S_i = \sum_{q \leq Q} \frac{q}{\phi(q)} \sum^*_{\chi (q)} \max_{y \leq x} |S_i|.$$
Easily we obtain
$$\sum_{q \leq Q} \frac{q}{\phi(q)} \sum^*_{\chi (q)} \max_{y \leq x} |\psi(y,\chi)| \leq \S_1+\S_2+\S_3+\S_4.$$
The heart of the proof of Theorem 1.3 in \cite{Akbary2015} are the following estimates:
\begin{lemma*}{\bf (Akbary, Hambrook, \cite[Section 7]{Akbary2015})}\label{akbary}
We have
\begin{equation*}
\begin{split}
& \S_1 \leq A_0 U Q^2, \;\; \S_2 <\left(x+Q^{\frac{5}{2}}V\right)(\log xV)^2,\;\; \S_3 < \S_3^{'}+\S_3^{''},\\
&\S_3^{'} < (x+Q^{\frac{5}{2}}U)(\log xU)^2,\\
&\S_3^{''} < \frac{c_3}{\log 2} \left( x + Qx^{\frac{1}{2}}U^{\frac{1}{2}}V^{\frac{1}{2}} + 2^{\frac{1}{2}}QxU^{-\frac{1}{2}}+Q^2x^{\frac{1}{2}}\right) (\log 2UV)^2(\log 4x),\\
&\S_4 < \frac{2^{\frac{3}{2}}A_1^{\frac{1}{2}c_3}}{\log 2} (x+QxV^{-\frac{1}{2}}+2^{\frac{1}{2}}QxU^{-\frac{1}{2}}+Q^2x^{\frac{1}{2}}) \left(\log \frac{2x}{V}\right)^{\frac{3}{2}}(\log e^3V)(\log 4x).
\end{split}
\end{equation*}
where $c_3$ as in (\ref{c3}).
\end{lemma*}
We estimate $\S_4$ contribution with the use of Lemma \ref{hh}. 
Writing $S_4$ as a dyadic sum we have
$$
S_4 = - \sum_{\bfrac{M = 2^{\alpha}}{\frac{1}{2}U < M \leq x/V}} \sum_{\bfrac{U < m \leq x/V}{ M < m \leq 2M}} \sum_{\bfrac{V < k \leq x/M}{ mk \leq y}} \Lambda(m) \left( \sum_{\bfrac{d \mid k}{ d \leq V}} \mu(d) \right) \chi(mk).
$$
Using the triangle inequality
$$
\S_4 \leq
\sum_{\bfrac{M = 2^{\alpha} }{ \frac{1}{2}U < M \leq x/V}}
\sum_{q \leq Q} \frac{q}{\phi(q)} \sum^*_{\chi (q)} \max_{y \leq x}
 \left| \sum_{\substack{U < m \leq x/V \\ M < m \leq 2M}}  
 \sum_{\substack{V < k \leq x/M \\ mk \leq y}} a_m b_k \chi(mk) \right|,
$$
where $a_m = \Lambda(m)$, and, as it was defined in the introduction $b_k = \sum_{{d | k},\;{d \leq V}} \mu(d)$.
Now apply the large sieve inequality (\ref{sieve}) to get
$$
\S_4 \leq c_3
\sum_{\bfrac{M = 2^{\alpha}}{\frac{1}{2}U < M \leq x/V }} 
(M'+Q^2)^{\frac{1}{2}}(K'+Q^2)^{\frac{1}{2}}\sigma_1(M)^{\frac{1}{2}} \sigma_2(M)^{\frac{1}{2}} L(M)
$$
where 
\begin{equation*}
\begin{split}
&\sigma_1(M) = \sum_{\bfrac{V < k \leq x/M }{}} |b_k|^{2},\;\;\;\;
\sigma_2(M)= \sum_{\bfrac{U < m \leq x/V }{ M < m \leq 2M}}|a_m|^2,\\
&\text{and} \\
&L(M) = \log \left( \frac{2x}{M} \min \left(\frac{x}{V},2M\right)  \right) \leq \log 4x,\\
\end{split}
\end{equation*}
where $M^{\prime}$ and $K^{\prime}$ denote the number of terms in the sums over $m$ and $k$, respectively. 
From the definition of $M'$ and $N'$ we conclude
\begin{equation*}
\begin{split}
& M'=\min\left(2M,\frac{x}{V}\right)-\max\left(M+1,U+1\right) \leq M,\\
& K'=\frac{x}{M}-(V+1)+1 \leq \frac{x}{M}.
\end{split}
\end{equation*}
By Chebyshev estimate we have an upper bound
$$\sigma_2(M) \leq \sum_{m \leq 2M} \Lambda(m)^2 \leq \psi(2M) \log 2M \leq 2A_0M\log 2M.$$
Thus by Cauchy inequality
\begin{equation} \label{thetas}
\begin{split}
&\S_4 \leq c_3 (\log 4x) \sum_{\bfrac{M = 2^{\alpha}}{\frac{1}{2}U < M \leq x/V }} (M+Q^2)^{\frac{1}{2}}\left(\frac{x}{M}+Q^2\right)^{\frac{1}{2}}(2A_0 M \log 2M)^{\frac{1}{2}} \sigma_1(M)^{\frac{1}{2}}. 
\end{split}
\end{equation}
Further
$$M(M+Q^2)\left(\frac{x}{M}+Q^2\right)=Mx+Q^2x+M^2Q^2+MQ^4$$
and
$$(\log 2M)^{\frac{1}{2}} \leq \left( \log \frac{2x}{V} \right)^{\frac{1}{2}}.$$
Using Lemma \ref{hh} we get
$$(\sigma_1(M))^{\frac{1}{2}} \leq \frac{x}{M}(L+C)-V(L+C)+2V^2,$$
that implies
$$\S_4 \leq c_3 (2A_0)^{\frac{1}{2}} (x+2^{\frac{1}{2}}Q^{\frac{1}{2}}xU^{-\frac{1}{2}}+QxV^{-\frac{1}{2}}+Q^2x^{\frac{1}{2}}) (\log 4x) \left( \log \frac{2x}{V} \right)^{\frac{1}{2}}\sum_{\bfrac{M = 2^{\alpha}}{\frac{1}{2}U < M \leq x/V }} 1.$$
Since 
$$\sum_{\bfrac{M = 2^{\alpha}}{\frac{1}{2}U < M \leq x/V }} 1 \leq \frac{\log \frac{2x}{V}}{\log 2},$$
then
$$\S_4 \leq \frac{c_3}{\log 2} (2A_0)^{\frac{1}{2}} (x+2^{\frac{1}{2}}Q^{\frac{1}{2}}xU^{-\frac{1}{2}}+QxV^{-\frac{1}{2}}+Q^2x^{\frac{1}{2}}) (\log 4x) \left( \log \frac{2x}{V} \right)^{\frac{3}{2}}.$$

Combining it with results of Lemma \ref{akbary} we get
\begin{equation} \label{polylog}
\S=\sum_{q \leq Q} \frac{q}{\phi(q)} {{\sum}^{*}}_{\chi (q)} \max_{y \leq x} |\psi(y,\chi)| \leq c_4 \Pol(x,Q,U,V) \Log(x,V,U),
\end{equation}
where
\begin{equation*}
\begin{split}
&c_4= \max \left\{ A_0,\frac{c_3}{\log 2},\frac{c_3}{\log 2} (2A_0)^{\frac{1}{2}}\right\} = \frac{c_3}{\log 2} (2A_0)^{\frac{1}{2}},\\
&\Pol(x,Q,U,V) = 4x + 2Q^2x^{\frac{1}{2}}+UQ^2+Q^{\frac{5}{2}}(U+V)+\\
&+2^{\frac{1}{2}}Q^{\frac{1}{2}}xU^{-\frac{1}{2}}+2^{\frac{1}{2}}QxU^{-\frac{1}{2}}+Qx^{\frac{1}{2}}U^{\frac{1}{2}}V^{\frac{1}{2}}+QxV^{-\frac{1}{2}},\\
&\Log(x,V,U) = \max \left\{ (\log xV)^2, (\log xU)^2, (\log 2UV)^2 \log 4x , \left( \log \frac{2x}{V}\right)^{\frac{3}{2}} \log 4x\right\},
\end{split}
\end{equation*}

Now let's specify $U$ and $V$. If $x^{\frac{1}{3}} \leq Q \leq x^{\frac{1}{2}}$, then $U = V = x^{\frac{2}{3}}{Q^{-1}}$.
Then putting that into previous expression we get for the factor
\begin{equation*}
\begin{split}
&\Pol_1(x,Q)= 4x+2Q^2x^{\frac{1}{2}}+Qx^{\frac{2}{3}}(1+2^{\frac{1}{2}})+Q^{\frac{3}{2}}x^{\frac{2}{3}}(2+2^{\frac{1}{2}}+1)+x^{\frac{7}{6}} \leq\\
&\leq 4x+2Q^2x^{\frac{1}{2}}+2Qx^{\frac{5}{6}}+Q^{\frac{3}{2}}x^{\frac{2}{3}}(2+2^{\frac{1}{2}}+1).
\end{split}
\end{equation*}
where we used the fact that $x^{\frac{7}{6}} \leq Qx^{\frac{5}{6}}$ and $Qx^{\frac{2}{3}} \leq Qx^{\frac{5}{6}}$.
Working in the same manner with $\Log$ and keeping in mind the condition $x \geq 4$ we find that
$$\Log_1(x,V,U)\leq \left(\frac{4}{3}\log x\right)^{\frac{3}{2}} 2 \log x = \frac{2^4}{3^{\frac{3}{2}}} (\log x)^{\frac{5}{2}}.$$
If $Q \leq x^{\frac{1}{3}}$, we let $U = V = x^{\frac{1}{3}}$ and get
\begin{equation*}
\begin{split}
&\Pol_2(x,Q)= 4x+2Q^2x^{\frac{1}{2}}+Q^2x^{\frac{1}{3}} + 2Q^{\frac{5}{2}} x^{\frac{1}{3}}+2^{\frac{1}{2}}Q^{\frac{1}{2}}x^{\frac{5}{6}}+Qx^{\frac{5}{6}}(2^{\frac{1}{2}}+2) \leq\\
&\leq x(5+2^{\frac{1}{2}})+2Q^2x^{\frac{1}{2}}+2Q^{\frac{3}{2}}x^{\frac{2}{3}} + Qx^{\frac{5}{6}}(2^{\frac{1}{2}}+2),
\end{split}
\end{equation*}
where we used $Q^2x^{\frac{1}{3}} \leq x$, $Q^{\frac{5}{2}}x^{\frac{1}{3}} \leq Q^{\frac{3}{2}}x^{\frac{2}{3}}$ and $Q^{\frac{1}{2}}x^{\frac{5}{6}} \leq x$.

Similarly we get for $$\Log_2(x,V,U) \leq 2 \left(\frac{7}{6}\right) ^{\frac{3}{2}} (\log x)^{\frac{5}{2}}.$$

Finally, we have in (\ref{polylog})
$$\S \leq c_4 \frac{2^4}{3^{\frac{3}{2}}} (7x+2Q^2x^{\frac{1}{2}}+5Q^{\frac{3}{2}} x^{\frac{2}{3}}+4Qx^{\frac{5}{6}})(\log x)^{\frac{5}{2}},$$
as demanded.

\section{Proof of Theorem 1}
Let $y \geq 2$,$(a,q)=1$.  By orthogonality of characters modulo $q$, we have
$$
\psi(y;q,a) = \frac{1}{\phi(q)} \sum_{\chi} \overline{\chi}(a) \psi(y,\chi).
$$
Define $\psi^{\prime}(y,\chi)=\psi(y,\chi)$ if $\chi \not = \chi_0$ and $\psi^{\prime}(y,\chi)= \psi(y,\chi) - \psi(y)$ otherwise,
$\chi_0$ is the principal character mod $q$. Then
\begin{align*}%\label{error decomposed}
\psi(y,q,a) - \frac{\psi(y)}{\phi(q)} = \frac{1}{\phi(q)} \sum_{\chi} \overline{\chi}(a) \psi^{\prime}(y,\chi).
\end{align*}
For a character $\chi$ (mod $q$), we let $\chi^{\ast}$ be the primitive character modulo $q^*$ inducing $\chi$.
Follow the way of \cite{Akbary2015} we obtain
$$
\psi^{\prime}(y,\chi^{\ast}) - \psi^{\prime}(y,\chi) 
= \psi(y,\chi^{\ast}) - \psi(y,\chi)
= \sum_{p^k \leq y} (\log p) (\chi^{\ast}(p^k) - \chi(p^k)).
$$
If $p | q$ then $(p^k,q^{*}) = 1$, and hence $\chi^{*}(p^k) = \chi(p^k)$.  If $p | q$ then $\chi(p^k) = 0$.  Therefore 
$$
|\psi^{\prime}(y, \chi^{\ast}) - \psi^{\prime}(y,\chi) |
\leq \sum_{\bfrac{p^k \leq y}{p | q}} (\log p) \leq (\log y) \sum_{p | q} 1 \leq (\log qy)^2.
$$
Denote the quantity we want to estimate as
$$\M=\sum_{\bfrac{q\leq Q }{ l(q)>Q_1}} \max_{2 \leq y \leq x} \max_{\bfrac{a}{(a,q)=1}} \left| \psi(y;q,a) - \frac{\psi(y)}{\phi(q)} \right|.$$
Since
$$
\left|\psi(y,q,a) - \frac{\psi(y)}{\phi(q)}\right|
\leq \frac{1}{\phi(q)} \sum_{\chi} |\psi^{\prime}(y,\chi)|
\leq (\log qy)^2 + \frac{1}{\phi(q)} \sum_{\chi} |\psi^{\prime}(y,\chi^{*})|,
$$
then
$$
\M \leq Q(\log Qx)^2 + \sum_{\bfrac{q \leq Q }{ l(q)>Q_1}}\frac{1}{\phi(q)} \sum_{\chi} \max_{2 \leq y \leq x} |\psi^{\prime}(y,\chi^{*})|.
$$
We have to take care just of the second term in the inequality above, since the first one is smaller than the desired bound.
It remains to prove 
$$\N=\sum_{\bfrac{q \leq Q }{ l(q)>Q_1}}\frac{1}{\phi(q)} \sum_{\chi} \max_{2 \leq y \leq x} |\psi^{\prime}(y,\chi^{*})| \leq (c_1-1)F(x,Q,Q_1)(\log x)^{4},$$
where $F(x,Q,Q_1)$ is the function from Theorem \ref{bvef}.
A primitive character $\chi^{*}\mmod q^{*}$ induces characters of moduli $dq^{*}$ and $\psi^{\prime}(y,\chi^{*}) = 0$ for $\chi$ principal, we observe
\begin{align*}
&\N = \sum_{\bfrac{q \leq Q }{l(q)>Q_1 }} \frac{1}{\phi(q)} \sum_{\bfrac{q^{*} | q }{ q^{*}\neq 1}} \sum^{*}_{\chi (q^{*})} \max_{2 \leq y \leq x}|\psi^{\prime}(y,\chi)| \leq \sum_{\bfrac{q^{*} \leq Q}{l(q^{*})>Q_1}} \sum^{*}_{\chi (q^{*})} \max_{2 \leq y \leq x} |\psi'(y,\chi)|
\sum_{k \leq \frac{Q}{q^{*}}}\frac{1}{\phi(kq^{*})}.\\
\end{align*}
As it was noted in \cite{Akbary2015} for $x>0$
$$ \sum_{k \leq x}\frac{1}{\phi(k)} \leq E_0 \log (ex)$$
and as $q^{*} \leq Q \leq x^{1/2}$, $\phi(k)\phi(q^{*}) \leq \phi(kq^{*})$ and $x \geq 4$, we have
$$
\sum_{k \leq \frac{Q}{q^{*}}} \frac{1}{\phi(kq^{*})} < \frac{5 E_0}{4 \phi(q^{*})} \log x.
$$
For $q>1$ and $\chi$ primitive character ($\mmod q$), we know that $\chi$ is non-principal and $\psi(y,\chi) = \psi'(y,\chi)$. Since we assumed $Q_1 \geq 1$ then we can can replace $\psi^{\prime}(y,\chi)$ by $\psi(y,\chi)$ inside the internal sum for $\N$.
Combining it with an expression for $\N$ we get
$$
\N \leq \frac{5 E_0}{4} (\log x) \sum_{\bfrac{q \leq Q}{l(q)>Q_1}} \frac{1}{\phi(q)}\sum_{\chi (q)}^{*} \max_{2 \leq y \leq x} |\psi(y,\chi)| = \Rc.
$$
Thus it remains to show that
$$\Rc \leq \frac{4(c_1-1)}{5E_0} F(x,Q,Q_1)(\log x)^{\frac{5}{2}}.$$
Let 
$$
\Rc(q) = \frac{q}{\phi(q)} \sum_{\chi(q)}^* \max_{2 \leq y \leq x} |\psi(y,\chi)|.
$$
Partial summation gives us
$$
\sum_{Q_1 < q \leq Q}\frac{1}{\phi(q)} \sum_{\chi (q)}^* \max_{2 \leq y \leq x} |\psi(y,\chi)| 
= \frac{1}{Q} \sum_{q \leq Q} \Rc(q) - \frac{1}{Q_1} \sum_{q \leq Q_1} \Rc(q) + \int_{Q_1}^{Q} \left(\sum_{q \leq t} \Rc(q)\right) \frac{dt}{t}. 
$$
Now we apply Theorem \ref{vaug}
$$\sum_{q\leq Q}\Rc(q) < c_0 f(x,Q)(\log x)^{\frac{5}{2}},$$
where $f(x,Q)=7x+2Q^2x^{\frac{1}{2}}+5Q^{\frac{3}{2}} x^{\frac{2}{3}}+4Qx^{\frac{5}{6}}$.
Then
$$
\sum_{Q_1 < q \leq Q}\frac{1}{\phi(q)} \sum_{\chi (q)}^* \max_{2 \leq y \leq x} |\psi(y,\chi)| 
< c_0\left(\Delta_f(Q,Q_1) + \int_{Q_1}^{Q}f(x,t)\frac{dt}{t}\right)(\log x)^{\frac{5}{2}}, 
$$
where $$\Delta_f(Q,Q_1) = \frac{f(x,Q)}{Q}-\frac{f(x,Q_1)}{Q_1} \leq \frac{7x}{Q_1}+2x^{\frac{1}{2}}Q+5x^{\frac{2}{3}}Q^{\frac{1}{2}}.$$
Calculating the integrals gives us
\begin{equation*}
\begin{split}
& \int_{Q_1}^Q f(x,t) \frac{dt}{t} < \frac{7x}{Q_1}+2x^{\frac{1}{2}}Q+10x^{\frac{2}{3}}Q^{\frac{1}{2}}+4x^{\frac{5}{6}}\log \frac{Q}{Q_1}.\\
\end{split}
\end{equation*}
Finally
$$
\N \leq \frac{4(c_1-1)}{5E_0} \left(\frac{14x}{Q_1}+4x^{\frac{1}{2}}Q+15x^{\frac{2}{3}}Q^{\frac{1}{2}}+4x^{\frac{5}{6}}\log \frac{Q}{Q_1}\right)(\log x)^{\frac{5}{2}}.
$$
\subsection{Proof of Remark \ref{pi}}
Define two functions
\begin{equation*}
\pi_{1}(y) = \sum_{2 \leq n \leq y} \frac{\Lambda(n)}{\log n} \;\;\;\;\;\;\; \text{and} \;\;\;\;\;\;\;
\pi_{1}(y;q,a) = \sum_{\bfrac{2 \leq n \leq y}{n \equiv a (\mmod q)}} \frac{\Lambda(n)}{\log n}.
\end{equation*}
Since
$$\pi_{1}(y;q,a)-\pi(y;q,a) = \sum_{2 \leq k \leq \frac{\log y}{\log 2}} \sum_{\bfrac{p^k \leq y}{p^k \equiv a (\mmod q)}} \frac{1}{k} \leq 
\sum_{2 \leq k \leq \frac{\log y}{\log 2}} \frac{\pi(y^{\frac{1}{2}})}{2} < 2 y^{\frac{1}{2}},$$
where we used the fact that for $x >1$(see for example \cite[Lemma 3.1]{Akbary2015})
$$\pi(x) < 1.25506 \frac{x}{\log x}.$$
Similarly,
$\pi_1(y) - \pi(y) < 2y^{\frac{1}{2}}. $
Thus by partial summation we obtain the bound
\begin{equation*}
\begin{split}
&\left|\pi_1(y;q,a) - \frac{\pi_1(y)}{\phi(q)} \right| = \left|\frac{\psi(y;q,a) - \psi(y) / \phi(q)}{\log y} -\int_{2}^{y} \frac{\psi(t;q,a) - \psi(t) / \phi(q)}{t \log^2 t} dt \right|  \\
&\leq \frac{1}{\log 2} \left|\psi(y;q,a) - \frac{\psi(y)}{\phi(q)}\right|+ \max_{2 \leq t \leq y} \left|\psi(t;q,a) - \frac{\psi(t)}{\phi(q)}\right| \left(\frac{1}{\log{2}}-\frac{1}{\log{y}}\right).
\end{split}
\end{equation*}
We have
\begin{equation*}
\begin{split}
&\sum_{\bfrac{q\leq Q }{ l(q)>Q_1}} \max_{2 \leq y \leq x} \max_{\bfrac{a}{(a,q)=1}} \left| \pi(y;q,a) - \frac{\pi(y)}{\phi(q)} \right| \\
&\leq \frac{2}{\log 2} \sum_{\bfrac{q \leq Q}{l(q)>Q_1}}  \max_{2 \leq y \leq x} \max_{a,(a,q)=1} \left| \psi(y,q,a) - \frac{\psi(y)}{\phi(q)} \right| + 2x^{\frac{1}{2}} \sum_{\bfrac{q \leq Q}{l(q)>Q_1}} \left(1+\frac{1}{\phi(q)}\right)\\
&<\frac{2c_1}{\log 2} F(x,Q,Q_1)(\log x)^{\frac{7}{2}} + 2x^{\frac{1}{2}} \sum_{\bfrac{q \leq Q}{l(q)>Q_1}} \left(1+\frac{1}{\phi(q)}\right),
\end{split}
\end{equation*}
where we used Theorem \ref{bvef} to estimate the first summand. For $x\geq 4$ 
$$2x^{\frac{1}{2}} \sum_{\bfrac{q \leq Q}{l(q)>Q_1}} \left(1+\frac{1}{\phi(q)}\right) < \frac{2c_1}{\log 2} F(x,Q,Q_1)(\log x)^{\frac{7}{2}}.$$
and we are done.

\bibliography{library}{}
\bibliographystyle{ieeetr}

\vspace{2em}

\medskip\noindent {\footnotesize Mathematisches Institut, Bunsenstraße 3-5, 37073, Göttingen, Germany
\hfil\break
e-mail: {\tt alisa.sedunova@phystech.edu}}

\end{document}